\documentclass[conference]{IEEEtran}
\IEEEoverridecommandlockouts
\usepackage{cite}
\usepackage{amsmath,amssymb,amsfonts}
\usepackage{algorithmic}
\usepackage{graphicx}
\usepackage{textcomp}
\usepackage{xcolor}
\def\BibTeX{{\rm B\kern-.05em{\sc i\kern-.025em b}\kern-.08em
    T\kern-.1667em\lower.7ex\hbox{E}\kern-.125emX}}
    
\newenvironment{rcases}
  {\left.\begin{aligned}}
  {\end{aligned}\right\rbrace}

\usepackage{url}
\usepackage{scalerel}

\newcommand{\argmin}{\operatornamewithlimits{argmin}}
\newcommand\abs[1]{\left|#1\right|}

\newcommand{\cm}[1]{\textcolor{black}{#1}}

\newcommand{\ve}[1]{\boldsymbol{#1}}

\begin{document}

\title{Community Energy Storage-based Energy Trading Management for Cost Benefits and Network Support \\}

\author{\IEEEauthorblockN{ Chathurika P. Mediwaththe}
\IEEEauthorblockA{\textit{College of Engineering and Computer Science} \\
\textit{The Australian National University}\\
Canberra, Australia \\
chathurika.mediwaththe@anu.edu.au}
\and
\IEEEauthorblockN{ Lachlan Blackhall}
\IEEEauthorblockA{\textit{College of Engineering and Computer Science} \\
\textit{The Australian National University}\\
Canberra, Australia \\
lachlan.blackhall@anu.edu.au}

}


\IEEEoverridecommandlockouts
\IEEEpubid{\makebox[\columnwidth]{978-1-7281-8550-7/20/\$31.00~\copyright2020 IEEE \hfill} \hspace{\columnsep}\makebox[\columnwidth]{ }}


\maketitle

\IEEEpubidadjcol 

\begin{abstract} 
In this paper, the extent to which the integration of rooftop photovoltaic (PV) power with a community energy storage (CES) system can reduce energy cost and distribution network (DN) loss is explored. To this end, three energy trading systems (ETSs) are compared; first, an ETS where PV users exchange energy with the CES system in addition to the grid, second, an ETS where PV users merely exchange energy with the CES system, and third, an ETS where PV users only exchange energy with the grid. A multi-objective optimization framework, combined with a linear distribution network power flow model, is developed to study the trade-off between the energy cost and network power loss reductions while satisfying the DN voltage and current flow limits. Simulations, with real energy demand and PV power data, highlight that enabling the energy exchange between the users and the CES system can give a better trade-off between the DN power loss and energy cost reductions. Further, simulations demonstrate that all three ETSs deliver nearly 85\% DN energy loss reduction with significantly increased revenues compared to an ETS without a CES system.

\end{abstract}


\section{Introduction}

The cost reductions of battery technologies have paved the way to integrate energy storage systems closer to users to facilitate effective energy management within residential communities. Community energy storage systems (CESs), in particular, are connected to the low voltage side of the distribution transformers in distribution networks (DNs) with the aim of supporting network operations. \cm{By exploiting the CES charge-discharge operation with the rooftop photovoltaic (PV) power, the network reliability can be increased by mitigating the technical challenges in low voltage DNs such as active power losses and voltage deviations beyond statutory range.}

Energy management problems with a CES system may involve competing objectives of different stakeholders, such as DN power loss reduction and voltage regulation for network operators and energy cost reduction for users. Hence, compared to the single objective optimization approach, multi-objective optimization formulation allows us to realize the trade-offs between several competing stakeholder objectives in an energy management framework with a CES system.

\cm{This paper studies the charge-discharge scheduling of a CES system with residential PV power by comparing three different energy trading systems (ETSs); (1) a system  where PV users can exchange energy with both the grid and the CES system (ETS 1), (2) a system where PV users can exchange energy only with the CES system (ETS 2), and (3) a system where PV users exchange energy only with the grid (ETS 3). By employing the linear version of the branch power flow model in \cite{Baran1}, a multi-objective optimization framework is developed to study the trade-off between the minimizations of the total energy costs for the users and the CES provider and DN power loss in the ETSs while satisfying the CES device limits, network current flow and voltage limits. }

In one branch of literature, optimization and control frameworks to exploit behind-the-meter user-owned energy storage systems to achieve economic and network benefits have been explored \cite{Appen,feeder_model, Hashemi, Bale, Celik, Cao}. Another branch of literature studies scheduling of centralized in-front-of-the-meter energy storage systems. To visualize the trade-off solutions among various network and economic benefits in the scheduling of centralized energy storage systems, multi-objective optimization frameworks have been proposed \cite{Tant, MG3, Karthikeyan2, Gao, Nick}. For instance, in \cite{Tant}, achieving network voltage regulation, peak demand reduction and the annual battery operating cost reduction has been studied by developing a multi-objective optimization framework. In \cite{Nick}, a multi-objective optimization framework has been explored to study the optimal trade-off between economic and technical goals including the DN voltage deviations, congestion, power losses, energy costs of supplying loads, and energy storage maintenance and investment costs. The common feature of the existing research is that the energy storage is merely enabled to exchange energy with the grid. In contrast to these single path energy flow-based frameworks, \cm{here, we explore the capability of sharing a CES system among multiple users by allowing them to trade energy directly with the storage system without breaching the network constraints such as voltage and current flow limits.}

\cm{We organize the paper as follows. The ETS configurations are given in Section~\ref{system_config} and the DN power flow model is given in Section~\ref{network_model}.  The multi-objective optimization framework is presented in Section~\ref{optimization_problem}. Numerical analyses and simulations are demonstrated in Section~\ref{results} with conclusions in Section~\ref{conclusion}.}

\section{Energy Trading System Models and Configurations}\label{system_config}
In the first two sub-sections of the this section, \cm{the models of the demand-side and the CES operation} used in the ETSs are generalized. Then, the three ETS configurations are explained according to how the energy flows between the CES system, users and the grid are enabled.

\subsection{Demand-side model}

\begin{figure}[t!]
\centering
\includegraphics[width=0.85\columnwidth]{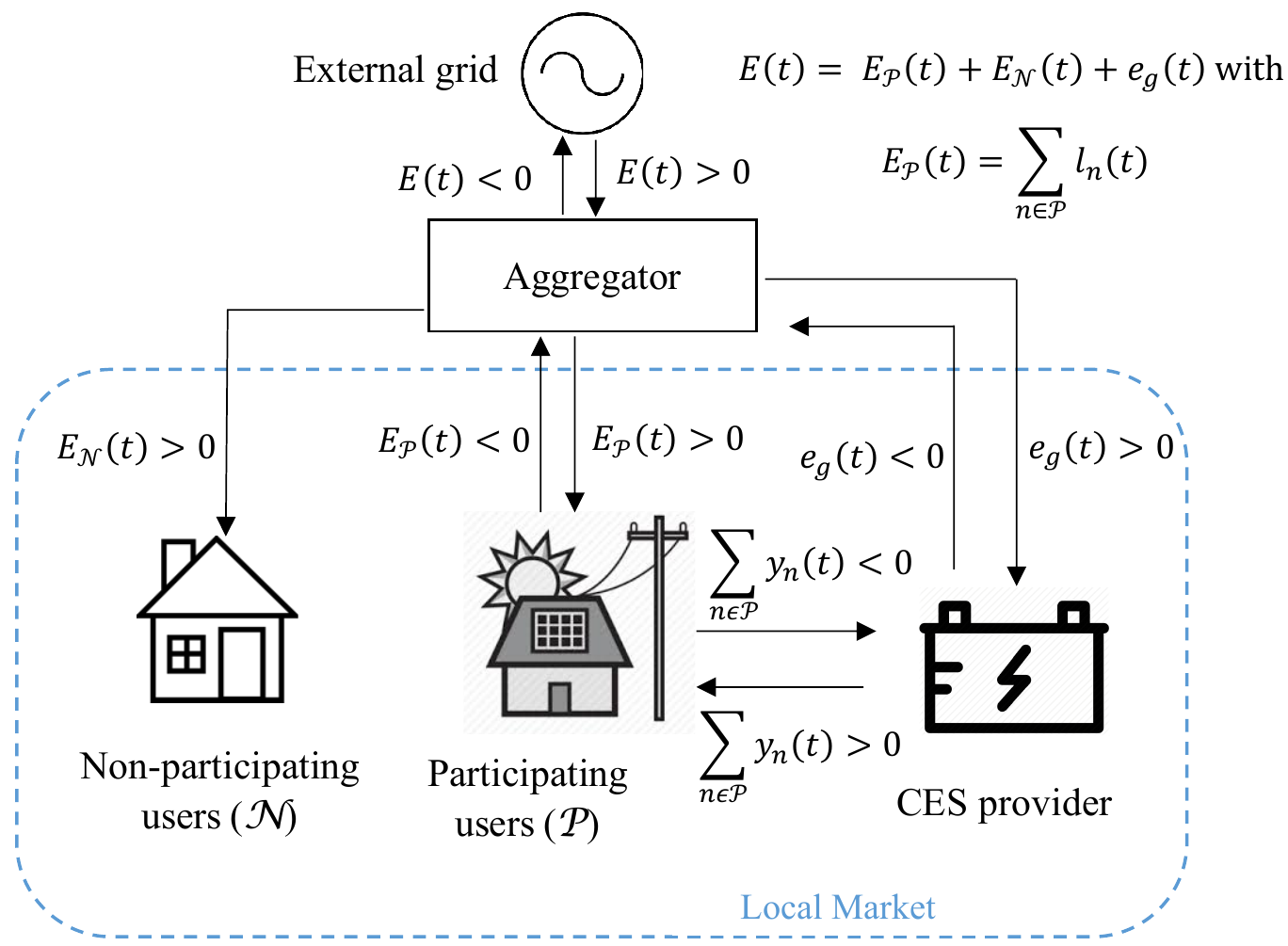}
\caption{ETS 1 - participating users exchange energy with both the CES device and the grid.}
\label{fig:sys1}
\end{figure}

\cm{Similar to our previous work in \cite{Power_sys_paper_chathurika}, an ETS model consists of a set of energy users $\mathcal{A}$, a CES system, and a system aggregator. We assume a third-party owns the CES system, and we refer to this owner as the CES provider. The energy users $\mathcal{A}$ are sub-divided into non-participating users $\mathcal{N}$ and participating users $\mathcal{P}$. Each user in $\mathcal{P}$ has rooftop solar and the users $\mathcal{N}$ do not have any power generation capability. We take PV systems at the users $\mathcal{P}$ operate at a power factor of 1. None of the users in $\mathcal{A}$ owns energy storage systems. The energy exchange between the local market and the grid is coordinated by the aggregator. Here, the local market, as shown in Fig.~\ref{fig:sys1}, comprises the CES system and the users $\mathcal{A}$.}

The time period of analysis $\mathcal{T}$ is split into $H$ number of time intervals. In this paper, $\mathcal{T}$ represents one day, and the time interval length is given by $\Delta t$. By considering demand and PV power generation variations, the users $\mathcal{P}$ are sub-divided into deficit users $\mathcal{P}^{-}(t)$ and surplus users $\mathcal{P}^{+}(t)$. At time $t$, $\mathcal{P} = \mathcal{P}^{-}(t)\cup \mathcal{P}^{+}(t)$. PV power at user $n\in \mathcal{P}$ at time $t$ is given by $g_n(t)\geq 0$, and  user $a\in \mathcal{A}$ has an energy demand $d_a(t)\geq 0$ at time $t$. Then, user  $n\in \mathcal{P}$ has an energy deficit $e_n(t)$ where $e_n(t) = d_n(t) - g_n(t)$. If $e_n(t)\geq 0$, then $n\in \mathcal{P}^{-}(t)$, and otherwise $n\in \mathcal{P}^{+}(t)$.


At time $t$, we consider user $n \in \mathcal{P}$ exchanges energy $y_n(t)$ with the CES and $l_n(t)$ with the grid. When buying, $y_n(t) >0$ and when selling, $y_n(t) <0$. The same sign convention is applied to $l_n(t)$. Then

\begin{equation}
e_n(t) = y_n(t) + l_n(t).\label{eq:id1}
\end{equation}

The aggregator determines the values for $l_n(t)$ day-ahead by solving a centralized optimization problem given in Section~\ref{optimization_problem} and using the reported $e_n(t)$ values by the users $\mathcal{P}$. Once $l_n(t)$ is known, $y_n(t)$ can be found by using the energy balance \eqref{eq:id1}. In this paper, it is supposed that the users $\mathcal{A}$ can obtain perfect information of energy demand and PV power generation forecasts for the next day.

\subsection{Community Energy Storage Model}
The CES system may charge-discharge with the grid in addition to charging and discharging with the users $\mathcal{P}$. The energy amount traded with the grid at time $t$ is denoted by $e_g(t)$. Here, $e_g(t) > 0$ represents buying energy and $e_g(t) < 0$ represents selling energy. It is supposed that $e_g(t),~\forall t \in \mathcal{T}$ are also determined by the aggregator by solving the optimization in Section~\ref{optimization_problem}. Given $y_n(t)$ and $e_g(t)$, $e_s(t) =  e_g(t) - \sum_{n\in \mathcal{P}}y_n(t)$ is the energy amount that is actually flowed in/out of the CES. Here, when charging, $e_s(t) > 0$, and when discharging $e_s(t) < 0$.

The constraint to satisfy the CES power ratings is given by
\begin{equation}
-\gamma^{\text{dis}}_{\text{max}} \leq \frac{e_s(t)}{\Delta t }\leq \gamma^{\text{ch}}_{\text{max}}, ~\forall{t \in \mathcal{T}} \label{eq:id2}
\end{equation}
where $\gamma^{\text{dis}}_{\text{max}}$ and $\gamma^{\text{ch}}_{\text{max}}$ are the maximum discharging and charging power rates of the CES, respectively. By incorporating the CES charging efficiency $0<\eta_c \leq 1$ and the discharging efficiency $\eta_d \geq 1$, $b(t) = b(t-1) + \eta e_s(t)$ gives the end of time $t$'s CES energy charge level. Here, $\eta = \eta_c$ if $e_s(t) \geq 0$, and $\eta = \eta_d$ otherwise. 
 For all $t\in \mathcal{T}$, the energy capacity constraint is 
\begin{equation}
B_{\text{min}}\leq b(t) \leq B_{\text{max}} \label{eq:id4}
\end{equation}
where $B_{\text{min}}$ and $B_{\text{max}}$ denote the minimum and maximum energy capacity limits of the CES, respectively.

We also take 
\begin{equation}
\abs{b(H) - b(0)} \leq  \theta \label{eq:id7}
\end{equation}
where $\theta$ is a small positive constant. This makes sure the continuity of the CES device's operation. 

\subsection{Different energy trading system configurations}\label{diff_config}
This section classifies the three different ETSs according to how the energy flows are enabled between the users $\mathcal{P}$, the grid, and the CES. 

\subsubsection{ETS 1} \label{all_paths}

In this system, the users $\mathcal{P}$ can trade energy with both the grid and the CES system. Hence, both $y_n(t)$ and $l_n(t)$ may have non-zero values. The ETS model is similar to the model in \cite{Power_sys_paper_chathurika}, however, in this paper, we consider a centralised optimisation approach to determine $l_n(t)$ whereas, in \cite{Power_sys_paper_chathurika}, the system utilises a distributed approach based on game theory.
By considering the amount of energy deficit at user $n\in \mathcal{P}$ at time $t$, the values of $l_n(t)$ are calculated such that
\begin{equation}
\begin{rcases}
               e_n(t) \leq  l_n(t) \leq  0,~~~\text{if}~n \in \mathcal{P}^{\scaleto{+}{4.5pt}}(t), \\
               0 \leq l_n(t) \leq e_n(t),~~~\text{if}~n \in \mathcal{P}^{\scaleto{-}{4.5pt}}(t).             
\end{rcases}\label{eq:id2a_b}
\end{equation}

\subsubsection{ETS 2}\label{WO_ln}

\begin{figure}[b!]
\centering
\includegraphics[width=0.85\columnwidth]{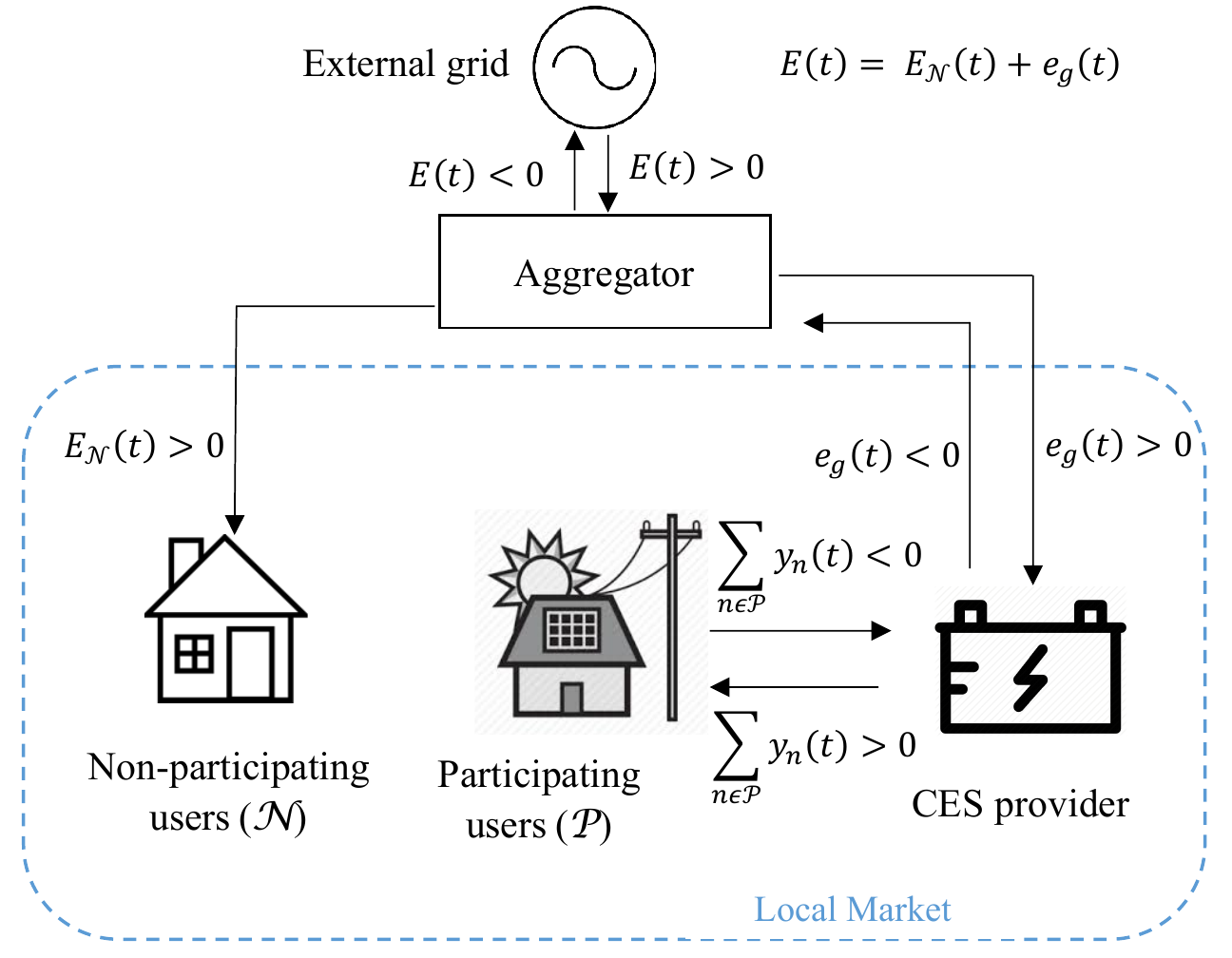}
\caption{ETS 2 - participating users exchanges energy only with the CES.}
\label{fig:sys2}
\end{figure}

In this system, participating users $\mathcal{P}$ interact only with the CES system and hence, exchange energy only with the storage. Therefore, only the energy transactions $y_n(t)$ exist, and 
\begin{equation}
l_n(t) = 0,~\forall n\in \mathcal{P},~\forall t \in \mathcal{T}. \label{eq:id2aa}
\end{equation}
Additionally, the CES system can exchange energy $e_g(t)$ with the grid (see Fig.~\ref{fig:sys2}).
Since the deficit energy of user $n\in \mathcal{P}$ is satisfied through the CES, $e_n(t) = y_n(t)$. 

\subsubsection{ETS 3}\label{WO_yn}

\begin{figure}[b!]
\centering
\includegraphics[width=0.85\columnwidth]{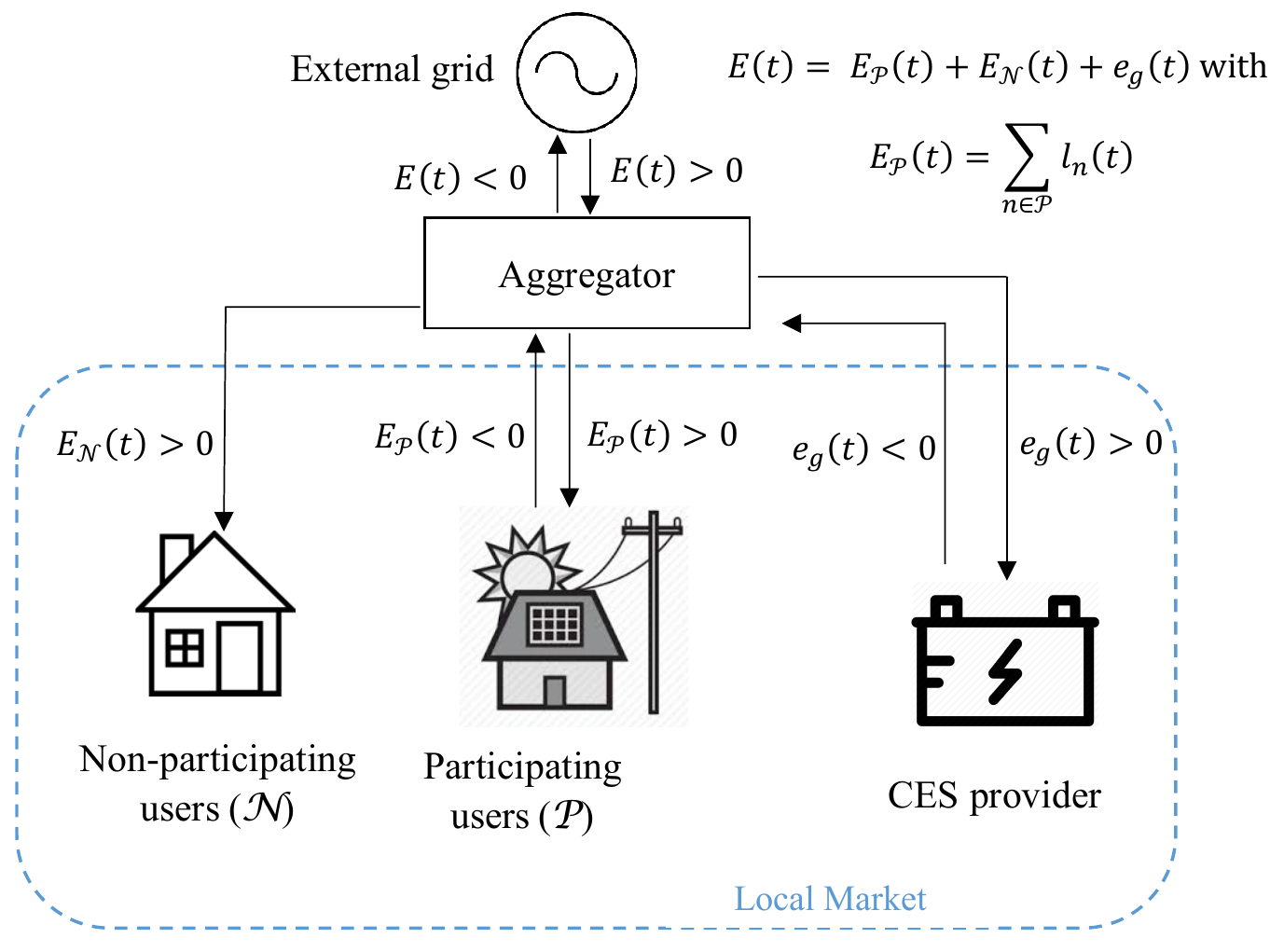}
\caption{ETS 3 - participating users exchanges energy only with the grid.}
\label{fig:sys3}
\end{figure}

In this system, the users $\mathcal{P}$ can only exchange energy with the grid and have no energy trading interaction with the CES system. Hence, only the energy transactions $l_n(t)$ exist, and $y_n(t) = 0,~\forall n\in \mathcal{P},~\forall t \in \mathcal{T}$. Since the deficit energy of user $n\in \mathcal{P}$ is satisfied through the grid, 
\begin{equation}
l_n(t) = e_n(t),~\forall n\in \mathcal{P},~\forall t \in \mathcal{T}. \label{eq:id2ab}
\end{equation}
Additionally, the CES can exchange energy $e_g(t)$ with the grid similar to the other two ETSs (see Fig.~\ref{fig:sys3}). Therefore, in this system, $e_s(t) =  e_g(t)$.

\section{Network Power flow Model}\label{network_model}
To calculate the DN voltages with the ETSs in Section~\ref{system_config}, this paper leverages the linearized Distflow equations in \cite{Baran1} developed for a radial DN.

Consider a radial DN described by a rooted tree $\mathcal{G}=(\mathcal{V}, \mathcal{E})$ where $\mathcal{V}= \{0,1,\dotsm,\text{N}\}$ is the set of buses and $\mathcal{E} = \{(i,j)\} \subset \mathcal{V}\times \mathcal{V}$ is the set of distribution lines. We ignore index $t$ in voltages and power in \eqref{eq:id22_a} for notational simplicity. 
Bus 0 is the slack bus representing the distribution transformer's secondary side. 
For bus $i \in \mathcal{V}$, $V_i$, $P_i$, and $Q_i$ are the voltage magnitude, active power consumption, and reactive power consumption, respectively. Additionally, $V_0$ is known and fixed. $r_{ij}$ and $x_{ij}$ are the resistance and reactance for each line $(i,j) \in \mathcal{E}$, respectively. Moreover, $P_{ij}$ and $Q_{ij}$ denote the active and reactive power flowing from bus $i$ to $j$. With $\mathbf{P} = (P_1,\dotsm,P_{\text{N}})^{\text{T}}$, $\mathbf{Q} = (Q_1,\dotsm,Q_{\text{N}})^{\text{T}}$, and $ \mathbf{V} = (V_1^2,\dotsm,V_{\text{N}}^2)^{\text{T}}$, the linearized Distflow model can be compactly written as \cite{Bitar}
\begin{equation}
 \mathbf{V} = -2\mathbf{R}\mathbf{P} -2 \mathbf{X}\mathbf{Q} + V_0^2\mathbf{1} \label{eq:id22_a}
\end{equation} 
Here, $\mathbf{R}, ~\mathbf{X}\in \mathbb{R}^{\text{N}\times \text{N}}$ where $\mathbf{R}_{ij} = \sum_{(h,k)\in \mathcal{L}_i \cap \mathcal{L}_j}r_{hk}$ and $\mathbf{X}_{ij} = \sum_{(h,k)\in \mathcal{L}_i \cap \mathcal{L}_j}x_{hk}$. $\mathcal{L}_i \subset \mathcal{E}$ is the set of lines on the path connecting bus $0$ and bus $i$. $\mathbf{1}$ is the N-dimensional all-ones vector.  
 
Given that $\mathbf{V}_{\text{min}}= V_{\text{min}}^2\mathbf{1}$ and $\mathbf{V}_{\text{max}} = V_{\text{max}}^2\mathbf{1}$ with $V_{\text{max}}$ and $V_{\text{min}}$ being the maximum and minimum DN voltage magnitude limits, respectively, it is required that 
\begin{equation}
\mathbf{V}_{\text{min}}~\leq~\mathbf{V}~\leq~\mathbf{V}_{\text{max}},~\forall t\in \mathcal{T}. \label{eq:id23_a}
\end{equation}
 
Without loss of generality, both non-participating and participating users may exist at a given bus $i \in \mathcal{V}\backslash \{0\}$. There the participating user set is given by $\mathcal{P}_i \subset \mathcal{P}$, and the non-participating user set is given by $\mathcal{N}_i \subset \mathcal{N}$. The active power consumption at bus $i$, $P_{i} = \frac{1}{\Delta t} \Big( \sum_{n\in \mathcal{P}_i}e_n(t) + \sum_{m\in \mathcal{N}_i}d_m(t)  + e_s(t)\Big)$ if the CES is placed at that bus. Otherwise, it is simply $P_{i} = \frac{1}{\Delta t} \Big( \sum_{n\in \mathcal{P}_i}e_n(t) + \sum_{m\in \mathcal{N}_i}d_m(t) \Big)$. Note that, here, $d_m(t)$ is the energy demand of user $m\in \mathcal{N}_i$. For each ETS, $e_s(t)$ is calculated accordingly with the choice of energy transactions $l_n(t)$, $y_n(t)$ and $e_g(t)$ as described in Section~\ref{diff_config}. If user $a\in \mathcal{A}_i$ has a reactive power demand $q_a(t)$, $Q_i$ can be taken as $Q_i= \sum_{a\in \mathcal{A}_i}q_a(t)$ where $ \mathcal{A}_i = \mathcal{P}_i \cup \mathcal{N}_i $. 

In addition to the voltage constraints, line current flow limits are also considered and given by
\begin{equation}
I^2_{ij}(t) \leq I_{ij,\text{max}}^2,~\forall t\in \mathcal{T},~\forall (i,j)\in \mathcal{E} \label{eq:id23_aa}
\end{equation}
where $I_{ij,\text{max}}$ is the magnitude of the maximum allowed complex current in line $(i,j)$, and $I^2_{ij}(t) = \frac{P_{ij}(t)^2 + Q_{ij}(t)^2 }{V^2_{0}(t)} $ by taking $V_{i}(t)^2 \approx V_{0}^2(t),~\forall i \in \mathcal{V} \backslash \{0\}$ \cite{Bitar}.

\section{Optimization Framework}\label{optimization_problem}

In each ETS, the aggregator solves a multi-objective optimization problem to find $l_n(t)$ and $e_g(t)$, and this section describes its generalized formation. The optimization includes two objectives: (1) the total energy cost, i.e., the sum of energy costs of the users $\mathcal{P}$ and the CES provider, minimization and (2) the DN power loss minimization.

By taking $\lambda_g(t)$ as the grid energy price, both the users $\mathcal{P}$ and the CES provider incur a total energy cost $f_{\text{cost}}$ that is given by
\begin{equation}
f_{\text{cost}} = \sum_{t=1}^H \lambda_g(t)\Big(E_{\mathcal{P}} (t)+ e_g(t) \Big). \label{eq:id9}
\end{equation}
where $E_{\mathcal{P}} (t) = \sum_{n\in \mathcal{P}}l_n(t)$.

Given $I^2_{ij}(t)$ as for \eqref{eq:id23_aa}, the total active power loss of the DN is calculated by using \cite{Bitar}
\begin{equation}
f_{\text{loss}} = \sum_{t=1}^H \sum_{(i,j) \in \mathcal{E}} r_{ij} I^2_{ij}(t). \label{eq:id24}
\end{equation}

After normalizing the objective functions and applying the linear-weighted combination \cite{mult_lin}, the multi-objective optimization problem is given by
\begin{equation}
 \ve{x^*} = \argmin_{\ve{x}~\in~\mathcal{X}}~\text{w}_1 \frac{f_{\text{cost}} - f_{\text{cost}}^{\text{utopia}}}{f^{\text{Nadir}}_{\text{cost}} - f_{\text{cost}}^{\text{utopia}}} + \text{w}_2 \frac{f_{\text{loss}}  - f^{\text{utopia}}_{\text{loss}}}{f^{\text{Nadir}}_{\text{loss}} - f^{\text{utopia}}_{\text{loss}}}\label{eq:id25}
\end{equation}
where $ \ve{x} = (\ve{E_{\mathcal{P}}}, \ve{e_g})$ with $\ve{E_{\mathcal{P}}} = (E_{\mathcal{P}}(1),\dotsm,E_{\mathcal{P}}(H))^{\text{T}}$ and $\ve{e_g} = (e_g(1),\dotsm,e_g(H))^{\text{T}}$. $\text{w}_1$ and $\text{w}_2$ are the weight coefficients, and $\mathcal{X}$ is the feasible set. The set $\mathcal{X}$ is subject to constraints \eqref{eq:id2a_b}, \eqref{eq:id1}-\eqref{eq:id7}, \eqref{eq:id23_a}-\eqref{eq:id23_aa} for the ETS 1, \eqref{eq:id2aa}, \eqref{eq:id1}-\eqref{eq:id7}, \eqref{eq:id23_a}-\eqref{eq:id23_aa} for the ETS 2, and \eqref{eq:id2ab}, \eqref{eq:id1}-\eqref{eq:id7}, \eqref{eq:id23_a}-\eqref{eq:id23_aa} for the ETS 3 in Section~\ref{diff_config}. Let the individual minimal points of $f_{\text{cost}}$ and $f_{\text{loss}}$ are given by $ \ve{x^*_{\text{cost}}} = \argmin _{{\ve{x}~\in~\mathcal{X}}} f_{\text{cost}}$ and $\ve{x^*_{\text{loss}}} = \argmin _{{\ve{x}~\in~\mathcal{X}}} f_{\text{loss}}$. Then the values for $f_{\text{cost}}^{\text{utopia}}$, $f_{\text{loss}}^{\text{utopia}}$, $f_{\text{cost}}^{\text{Nadir}}$, and $f_{\text{loss}}^{\text{Nadir}}$ in \eqref{eq:id25} are calculated as follows;
\begin{equation*}
\begin{rcases}
               f_{\text{cost}}^{\text{utopia}} = f_{\text{cost}}(\ve{x^*_{\text{cost}}}), \\
              f_{\text{loss}}^{\text{utopia}}= f_{\text{loss}}(\ve{x^*_{\text{loss}}}), \\
              f_{\text{cost}}^{\text{Nadir}} = \text{max} [f_{\text{cost}}(\ve{x^*_{\text{cost}}})~f_{\text{cost}}(\ve{x^*_{\text{loss}}})],\\
              f_{\text{loss}}^{\text{Nadir}} = \text{max} [f_{\text{loss}}(\ve{x^*_{\text{cost}}})~f_{\text{loss}}(\ve{x^*_{\text{loss}}})].
\end{rcases}\label{eq:id2a}
\end{equation*}
To find the coefficients $\text{w}_1$ and
$\text{w}_2$, the Analytic Hierarchy Process
(AHP) \cite{Saaty2004} can be applied, and $\sum_{k=1}^2 \text{w}_k =1$ with $\text{w}_k \in [0,1]$. If the set $\mathcal{P}$ comprises either only deficit users or surplus users, we take $l_n^*(t) = E_{\mathcal{P}}^*(t)\frac{e_n(t)}{\sum_{\mathcal{P}}e_n(t)}$ once $E_{\mathcal{P}}^*(t)$ is found by solving \eqref{eq:id25}. If both types of users exist in $\mathcal{P}$, then we set $l_n^*(t) = 0$.
\eqref{eq:id25} is a convex quadratically-constrained quadratic program, and thus, \eqref{eq:id25} can be solved by using convex optimization algorithms in \cite{boyd2004convex}. 

\section{Numerical analyses and simulations}\label{results}
\begin{figure}[b!]
\centering
\includegraphics[width=0.85\columnwidth]{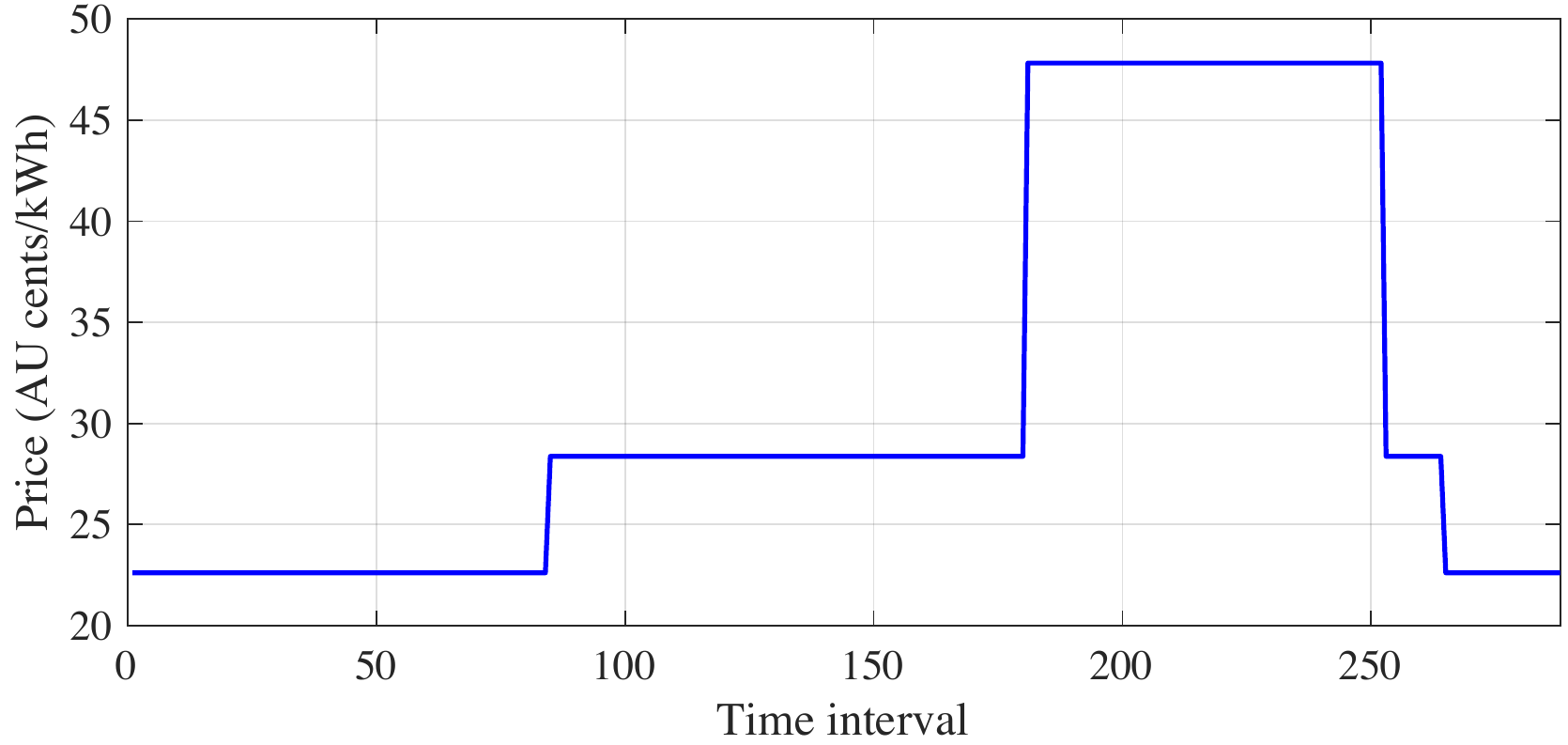}
\caption{Grid energy price signal $\lambda_g(t)$.}
\label{fig:price}
\end{figure}

In simulations, the radial distribution feeder with 7 buses given in \cite{feeder_model} is considered.
\begin{figure}[t!]
\centering
\includegraphics[width=0.85\columnwidth]{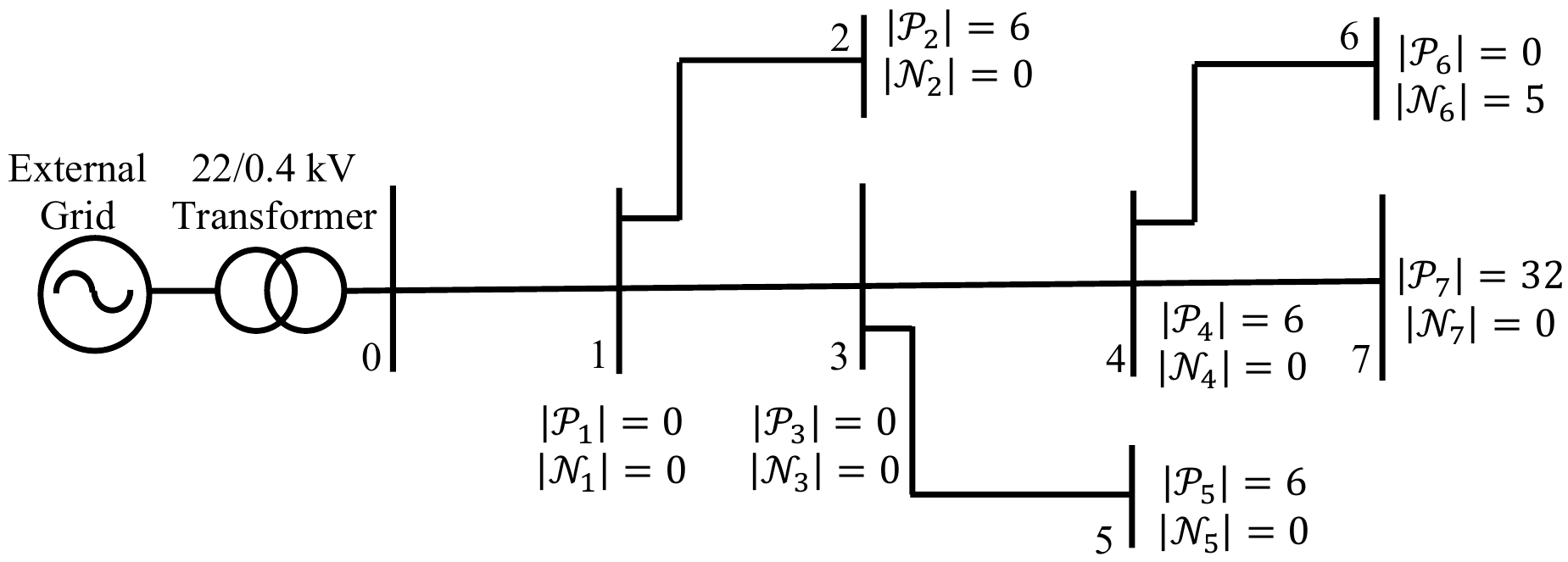}
\caption{7-bus radial feeder with the number of users at each bus.}
\label{fig:LV_feeder}
\end{figure}
55 users are connected to the feeder as shown in Fig.~\ref{fig:LV_feeder}, and hence, $\abs{\mathcal{A}} = 55$. Real PV power and active power demand data of 55 users in Canberra, Australia \cite{NextGen} are used to generate users' daily active power and PV power profiles in simulations. The generated daily power profiles represent the selected users' average Autumn daily demand and PV power profiles. Reactive power demand is not considered because of lack of real data. Additionally, $B_{\text{max}}= 900~\text{kWh},~B_{\text{min}}=0.05B_{\text{max}},~\gamma^{\text{ch}}_{\text{max}}=\gamma^{\text{dis}}_{\text{max}}=400~ \text{kW},~\eta_d = 1.02,~ \eta_c = 0.98,~\Delta t = 1/12~\text{hrs},~H = 288,~V_{\text{0}} = 1~\text{p.u.},~V_{\text{min}} = 0.95~\text{p.u.},~V_{\text{max}} = 1.05~\text{p.u.}$. The CES is placed at bus 7. $\text{w}_1 = 0.67$ and $\text{w}_2 = 0.33$ which are found by using the pairwise comparison matrix in \cite{Nick} with the AHP. The grid energy price signal $\lambda_g(t)$ is taken as the the time-of-use price signal in \cite{Origin} and is shown in Fig.~\ref{fig:price}.
In this paper, to compare the performance of the three ETSs, a baseline is considered without a CES. The users $\mathcal{P}$ in the baseline exchanges energy only with the grid for which they receive a price $\lambda_g(t)$.

\begin{figure}[b!]
\centering
\includegraphics[width=0.85\columnwidth]{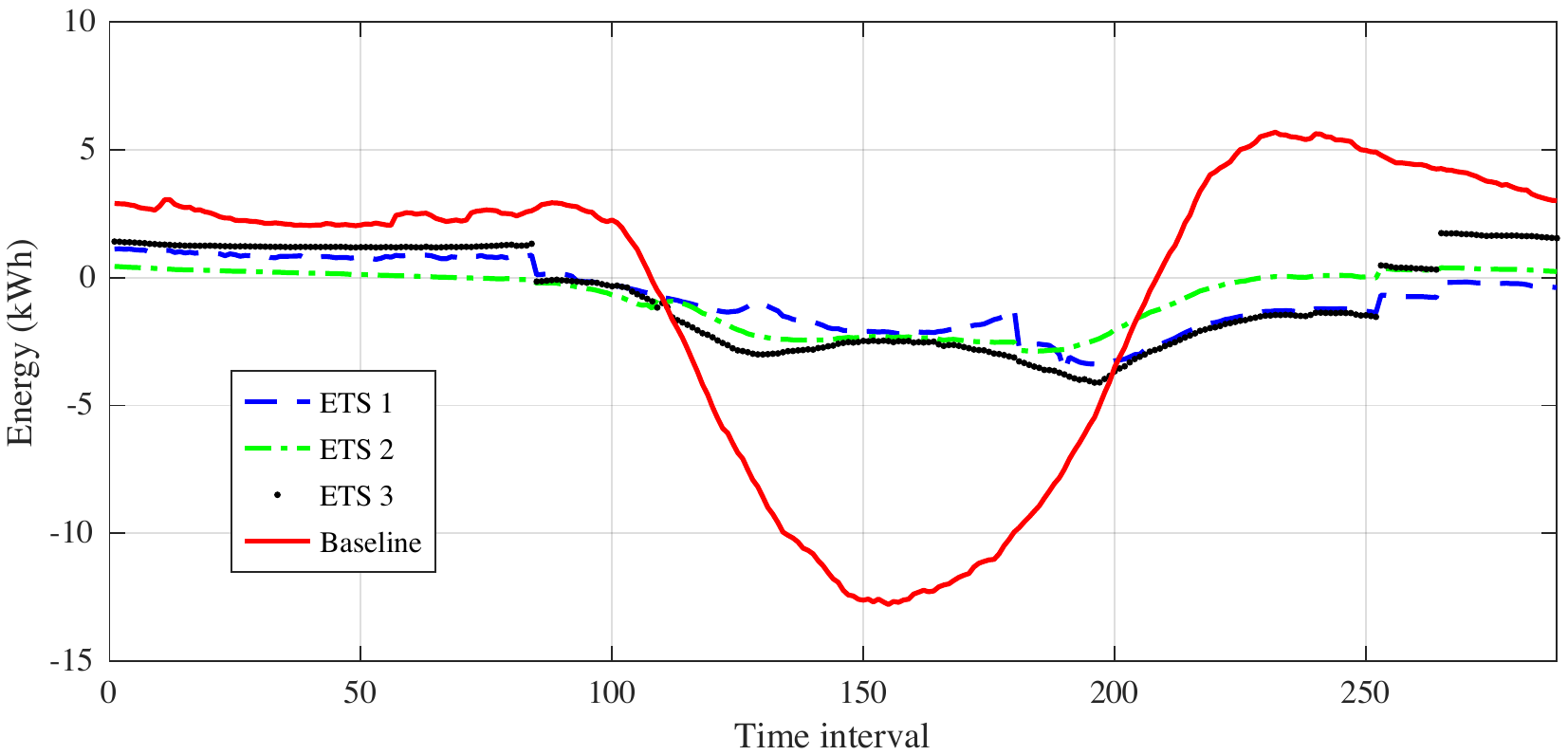}
\caption{Temporal variations of the total grid energy load $E(t)$.}
\label{fig:Grid energy}
\end{figure}

\begin{figure}[t!]
\centering
\includegraphics[width=0.85\columnwidth]{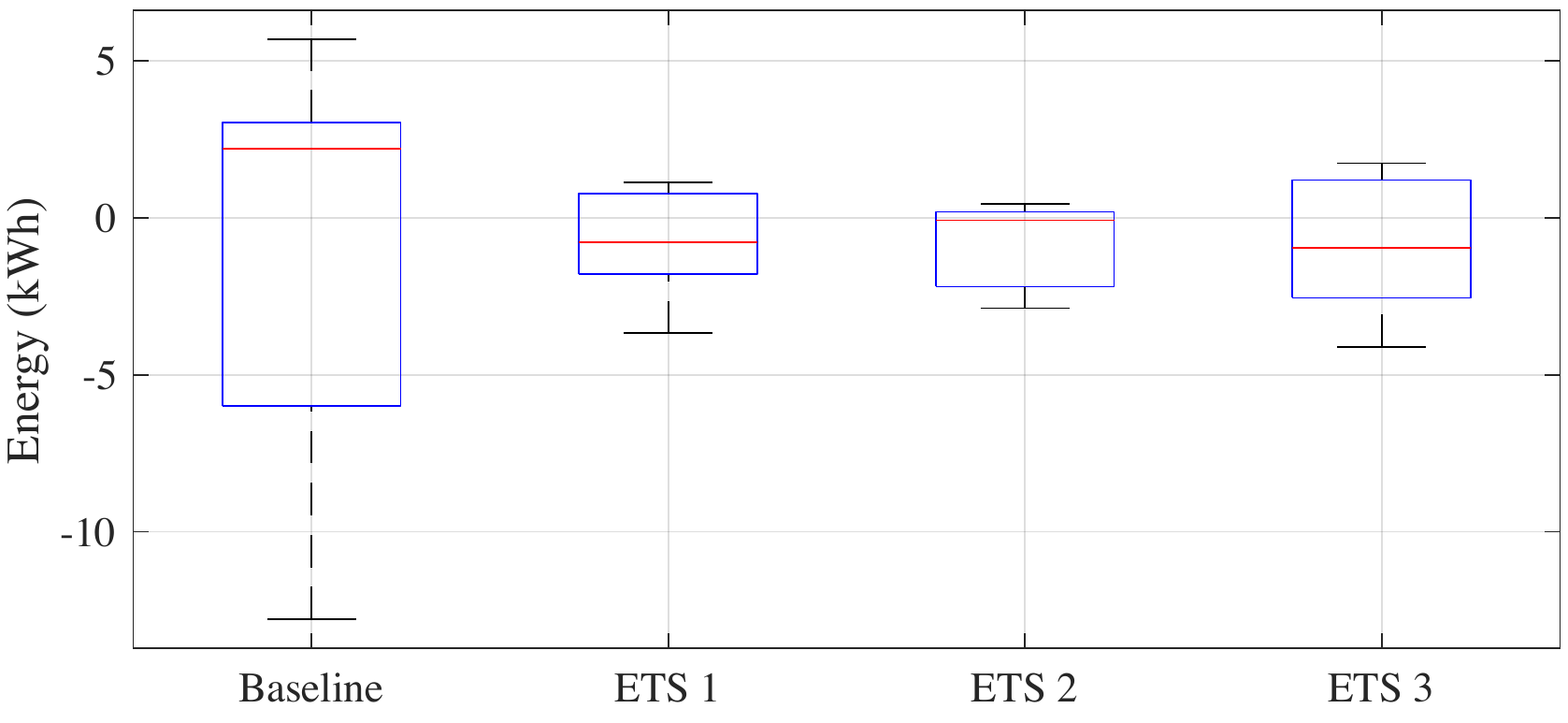}
\caption{Distributions of 24-hr total grid energy load $E(t)$.}
\label{fig:boxplot grid energy}
\end{figure}

Fig.~\ref{fig:Grid energy} illustrates the temporal variations of the total grid energy load $E(t)$, that is given by $E(t) = E_{\mathcal{P}}(t) + E_{\mathcal{N}}(t) + e_g(t)$ where $E_{\mathcal{N}}(t)$ is the total grid energy of the users $\mathcal{N}$ at time $t$, in the three ETSs and in the baseline system. 
As shown in the figure, during the day, $E(t)$ is negative in the baseline system due to excess PV power generation at the users $\mathcal{P}$.  
Fig.~\ref{fig:boxplot grid energy} depicts the distributions of the total grid load $E(t)$ within the 24-hr time period in the three ETSs and in the baseline system. 
In this figure, the maximum and the minimum grid loads are illustrated by the whisker horizontal bars. The figure depicts the ETSs reduce the positive peak grid energy demand and the peak reverse power flow to the grid (negative grid load) that occurs due to excess PV power during the day in the baseline. Out of the three ETSs, the ETS 2 results in the lowest peak positive and negative grid energy loads as the users $\mathcal{P}$ predominantly exchange all their energy transactions with the CES system.  

\begin{figure}[b!]
\centering
\includegraphics[width=0.85\columnwidth]{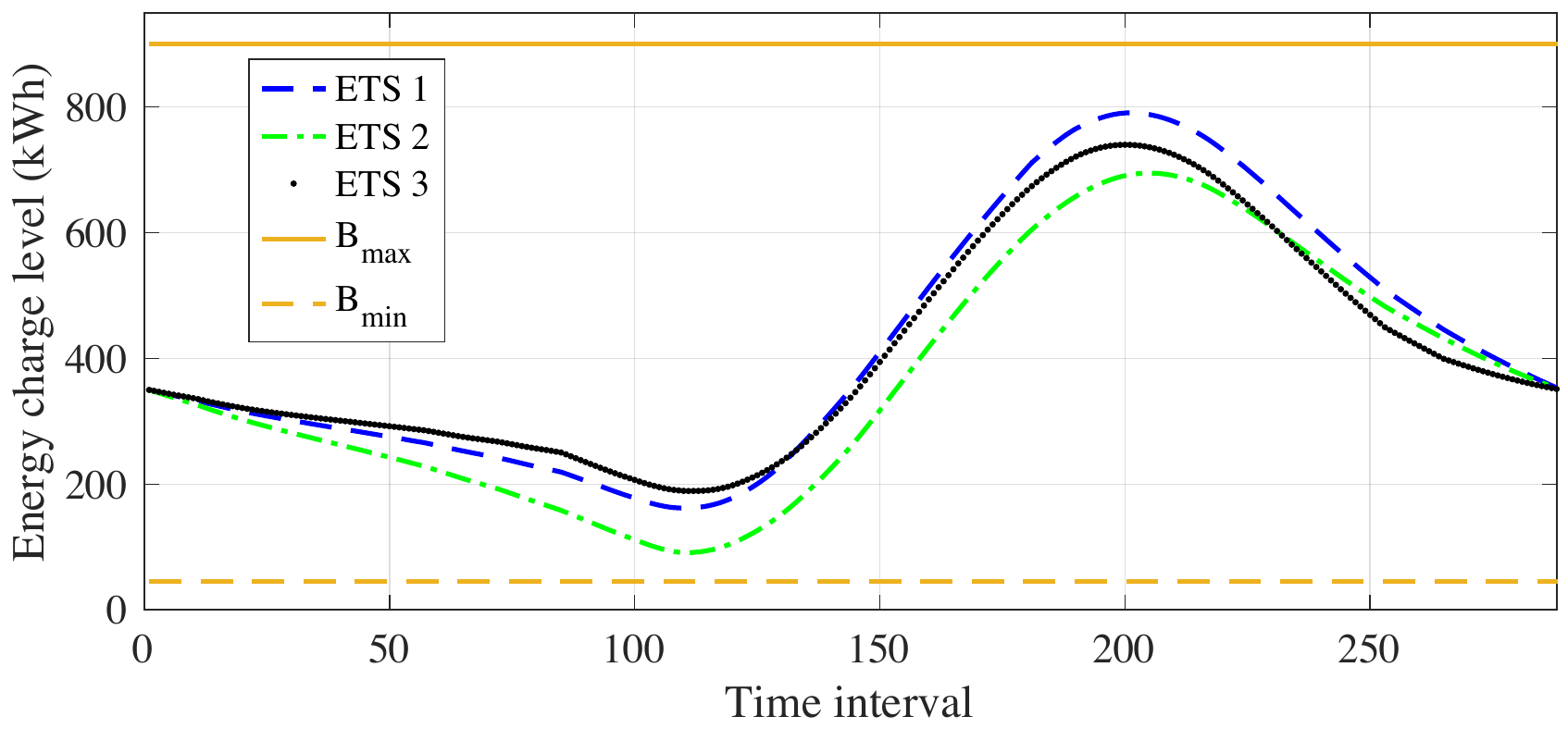}
\caption{The variations of the CES system's energy charge levels.}
\label{fig:CES SOC}
\end{figure}

\begin{figure}[t!]
\centering
\includegraphics[width=0.85\columnwidth]{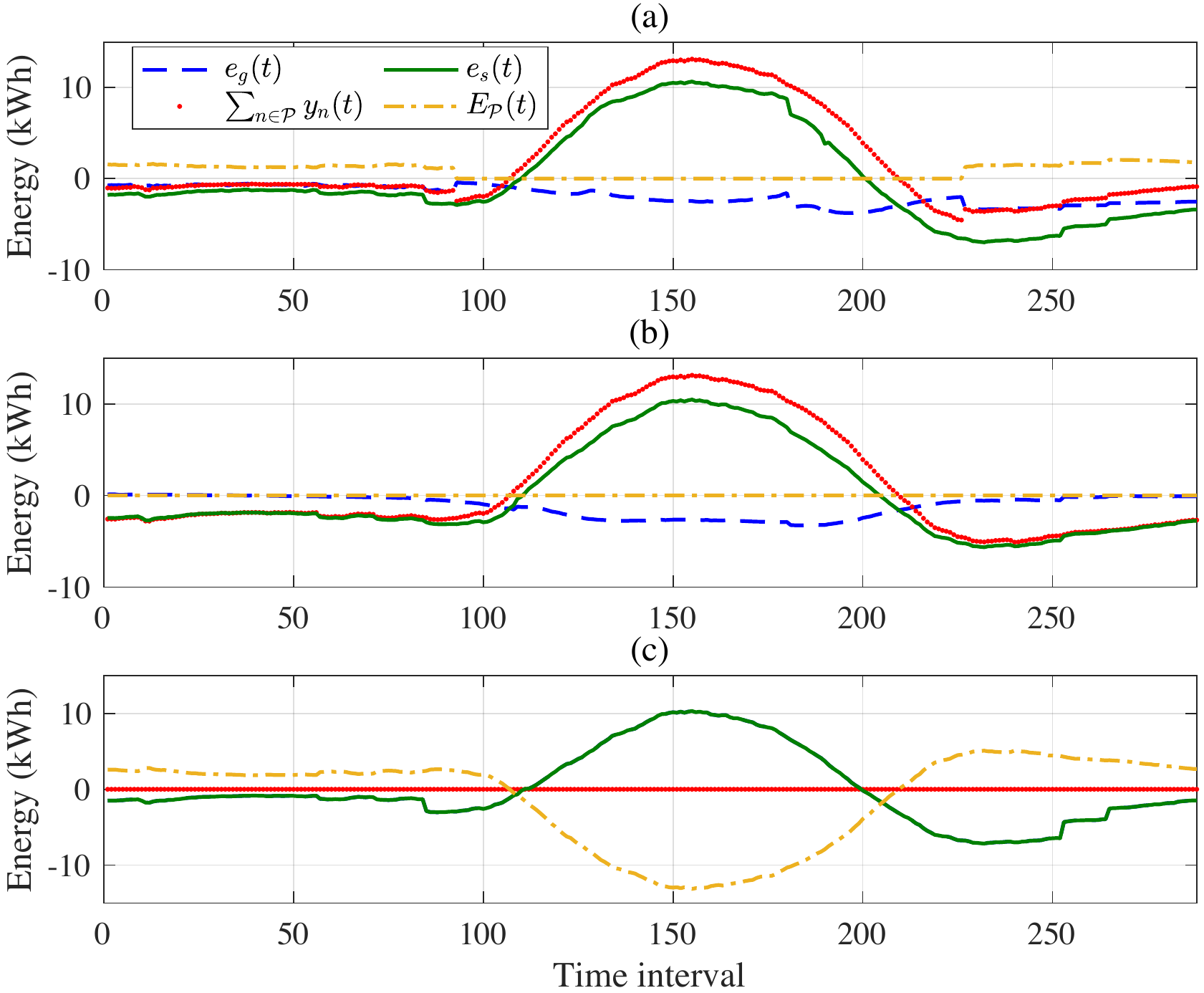}
\caption{Energy transactions of the (a) ETS 1 (b) ETS 2 and (c) ETS 3.}
\label{fig:energy transactions all paths}
\end{figure}

Fig.~\ref{fig:CES SOC} shows the temporal variations of the CES energy charge levels with the three ETSs. 
Additionally, Figs.~\ref{fig:energy transactions all paths}(a), \ref{fig:energy transactions all paths}(b) and \ref{fig:energy transactions all paths}(c) depict the different energy transactions in the three ETSs.
Fig.~\ref{fig:CES SOC} shows, in all systems, the CES system is charging during the day ($e_s(t) > 0$) and discharging in the afternoon ($e_s(t) < 0$). 
As shown in Figs. \ref{fig:energy transactions all paths}(a) and \ref{fig:energy transactions all paths}(b), in the ETSs 1 and 2, the CES system exchanges energy with both the users $\mathcal{P}$ and the grid, i.e., the profiles of $\sum_{n\in \mathcal{P}}y_n(t)$ and $e_g(t)$ have non-zero values. In these two systems, the CES charges during the day mainly due to the PV energy sold by the users ($\sum_{n\in \mathcal{P}}y_n(t) >0$). However, in ETS 3, because the users do not exchange energy with the CES, i.e., $\sum_{n\in \mathcal{P}}y_n(t) =0$,  the CES charge-discharge profile $e_s(t)$ coincides with the CES provider's grid energy trading profile $e_g(t)$ as shown in  Fig.~\ref{fig:energy transactions all paths}(c). Therefore, in the ETS 3, the charging of the CES during the day occurs mainly due to the energy bought by the CES provider ($e_g(t) > 0$). In all three systems, the CES system discharges in the afternoon and that reduces the afternoon peak grid energy demand as shown in Fig.~\ref{fig:Grid energy}. 

\begin{figure}[t!]
\centering
\includegraphics[width=0.85\columnwidth]{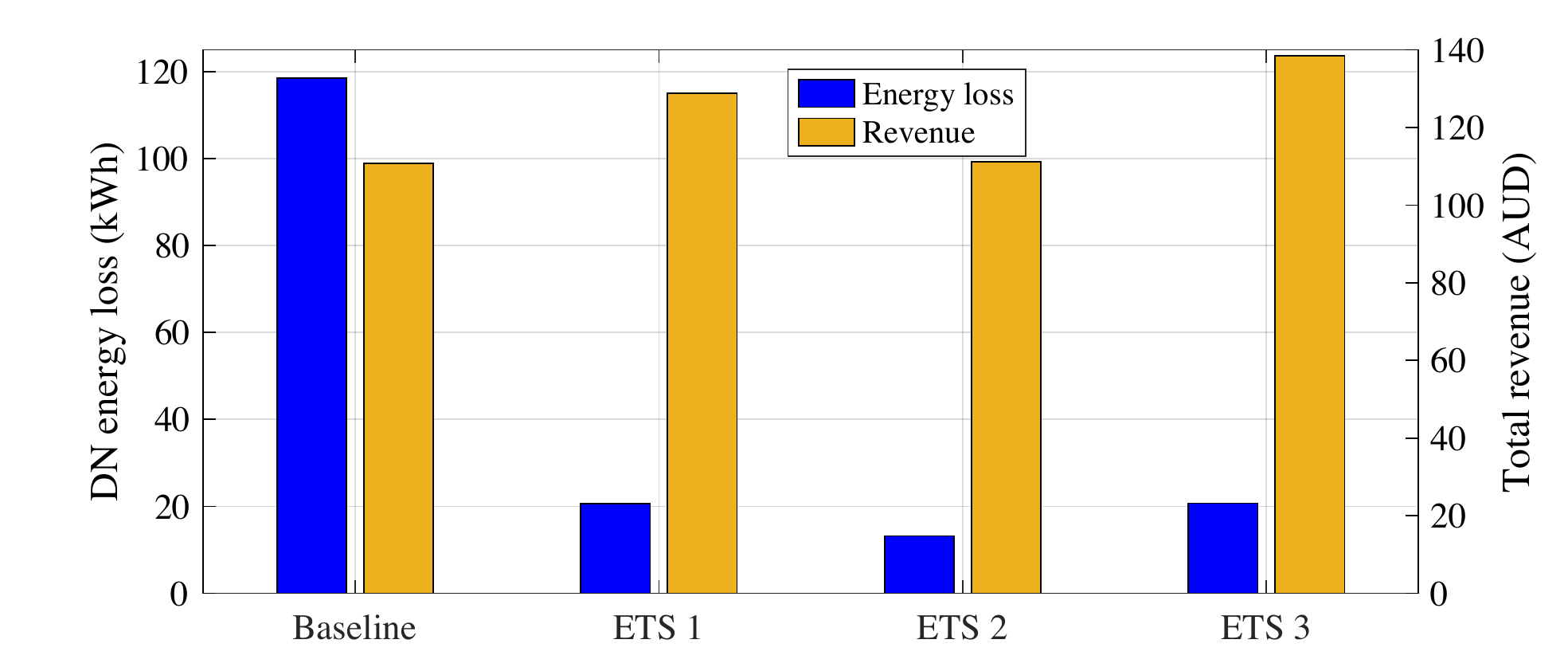}
\caption{DN losses and total revenues of the ETSs and the baseline.}
\label{fig:bar graph DN loss and revenue}
\end{figure}

Fig.~\ref{fig:bar graph DN loss and revenue} compares the total revenue ($-f_{\text{cost}}$ in \eqref{eq:id9}) and the DN energy loss, by using \eqref{eq:id24}, in the three ETSs with those of the baseline system. 
As shown, all three ETSs are capable of reducing the DN energy loss of the baseline system by nearly 85\%. Additionally, the total revenues of the three ETSs are higher than that in the baseline system. The total revenue decreases from the ETS 3 through the ETS 1 to the ETS 2. A similar trend was observed for the DN energy loss.  In the ETS 3, both participating users and the CES provider receive positive cumulative revenues and that leads to the greatest revenue of the three systems. On the other hand, in the ETS 2, the users receive zero revenue because they trade energy only with the CES system ($l_n(t) = 0$), without receiving a price for their energy transactions. That leads to the least revenue in the ETS 2. However, the reduction of the revenue of the ETS 2 compared to the revenue of the ETS 3, where the users only exchange energy with the grid, is only by 19\%, and the reduction of the revenue of the ETS 1, where the users trade energy with both the grid and the CES system, compared to revenue of the ETS 3 is only by 7\%. These trends highlight that energy trading between the users and the CES system can give a better trade-off between the DN power loss and energy cost reductions.

\begin{figure}[b!]
\centering
\includegraphics[width=0.85\columnwidth]{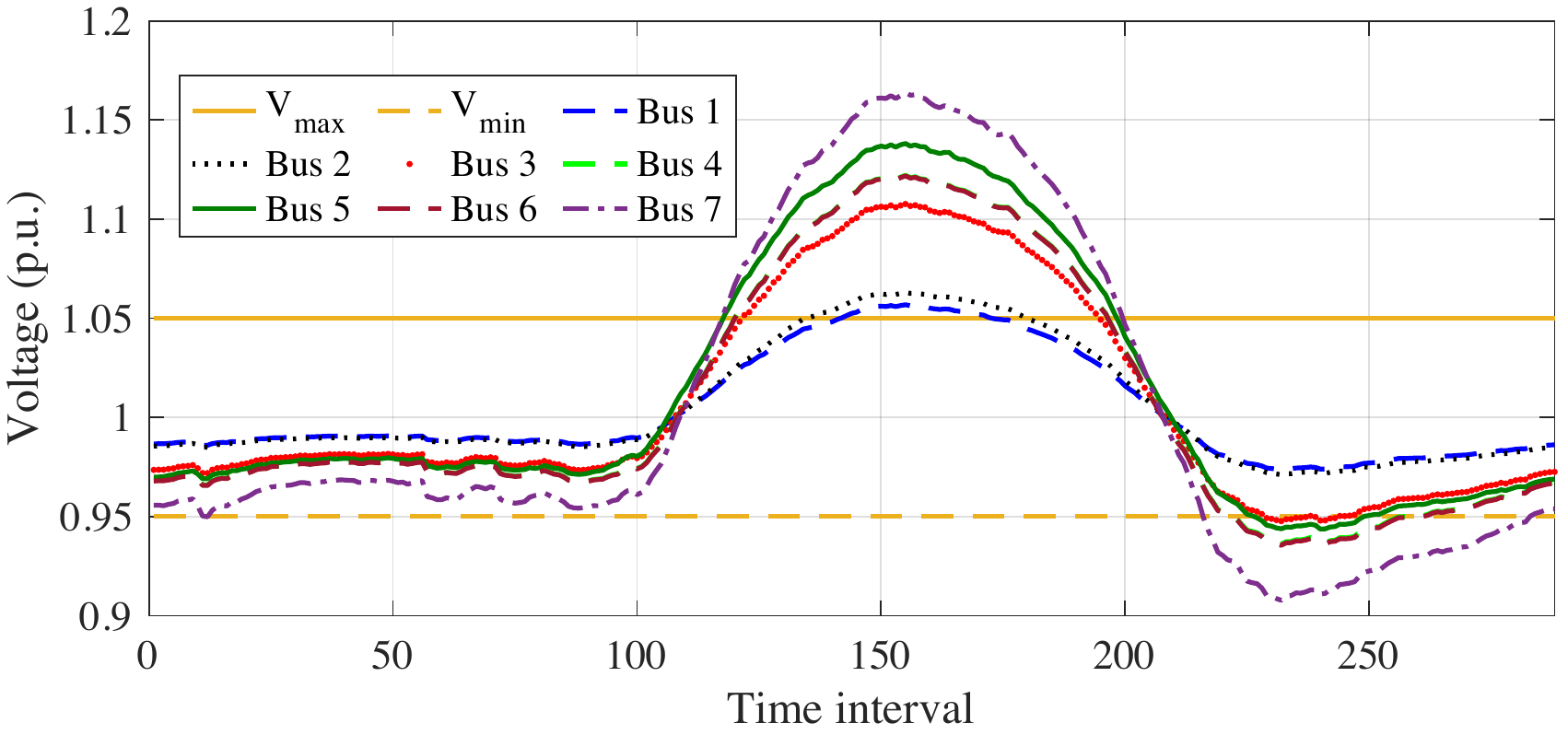}
\caption{Bus voltage profiles of the baseline system.}
\label{fig:voltages baseline}
\end{figure}

\begin{figure}[t!]
\centering
\includegraphics[width=0.85\columnwidth]{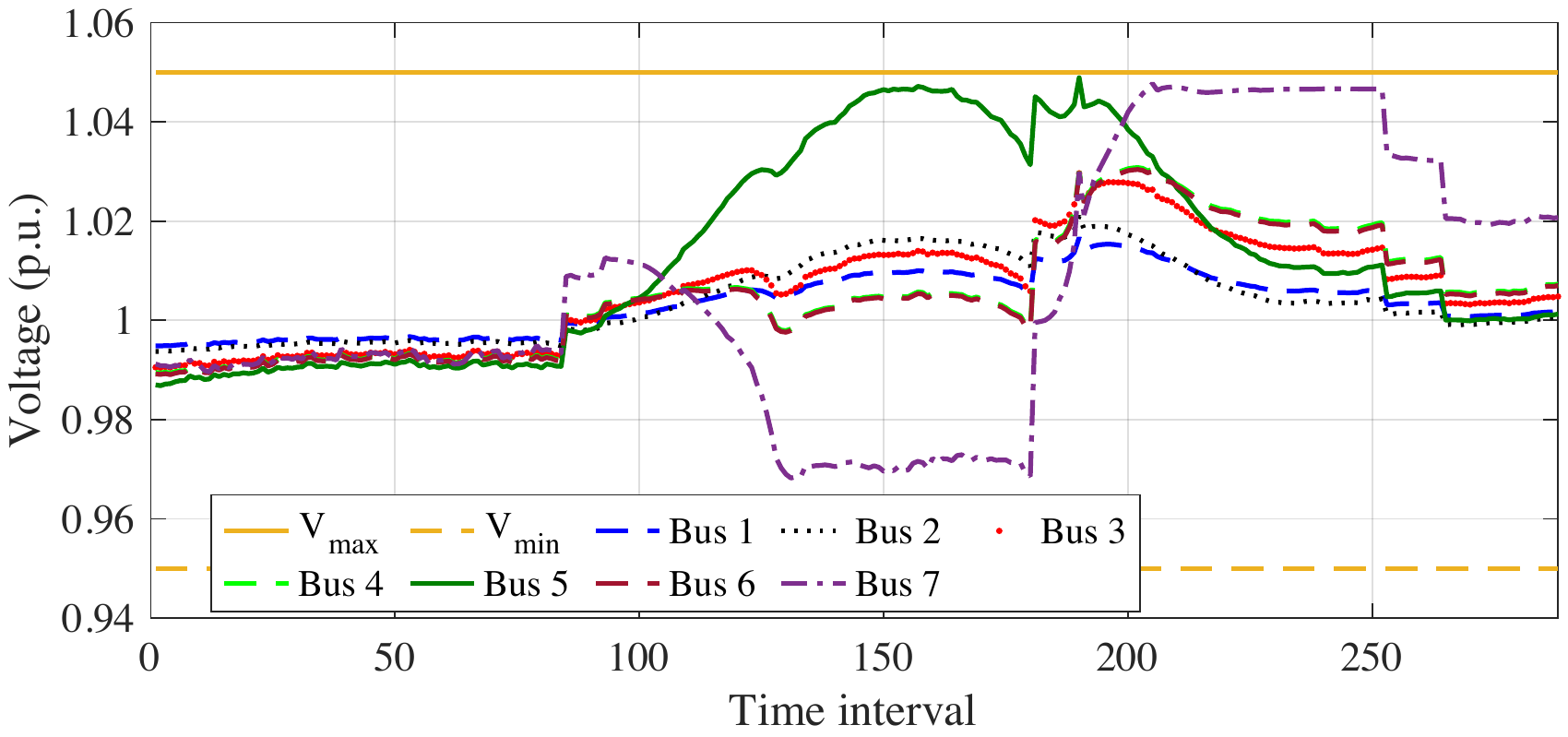}
\caption{Bus voltage profiles of the ETS 1.}
\label{fig:voltages all paths}
\end{figure}

\begin{figure}[t!]
\centering
\includegraphics[width=0.85\columnwidth]{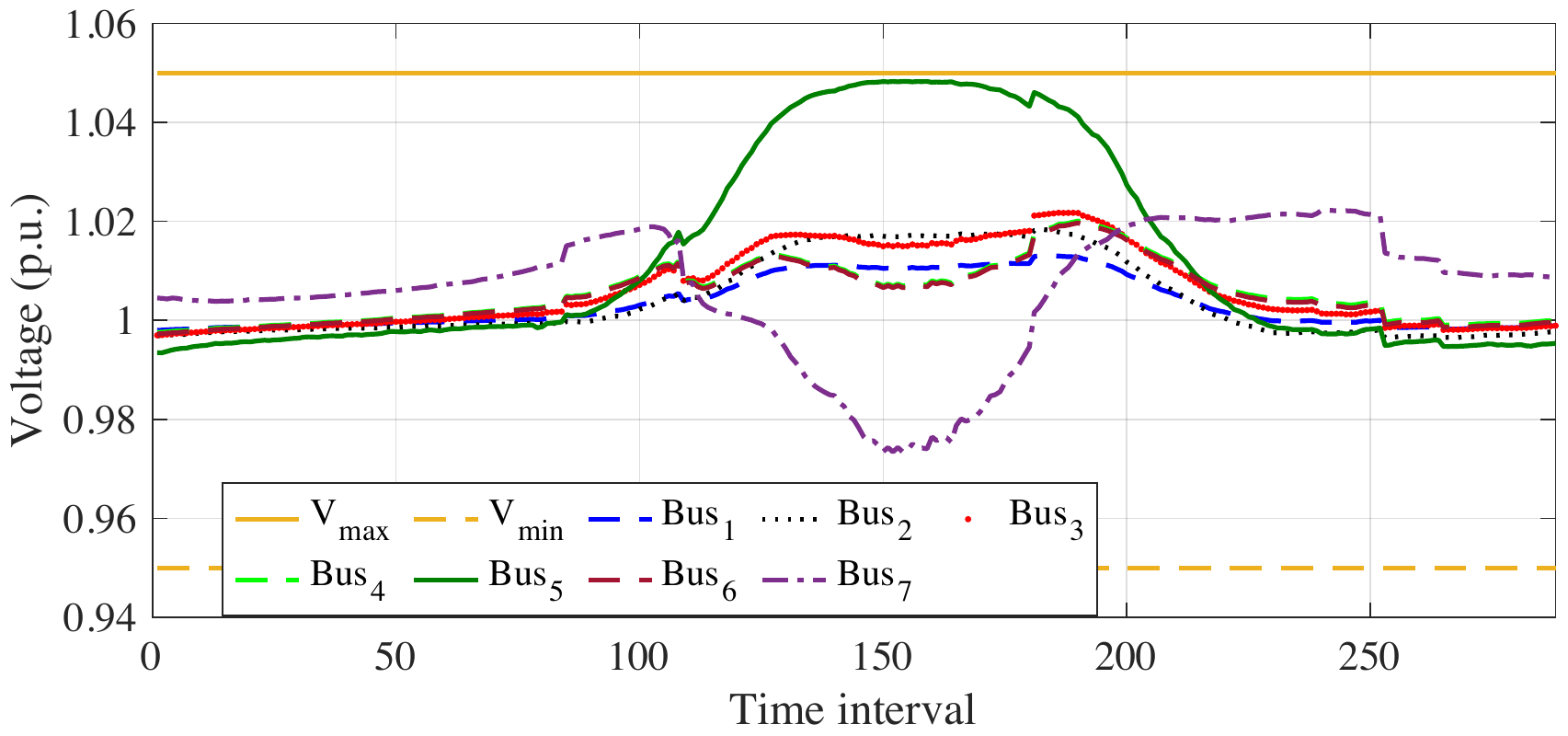}
\caption{Bus voltage profiles of the ETS 2.}
\label{fig:voltages without ln(t)}
\end{figure}

\begin{figure}[t!]
\centering
\includegraphics[width=0.85\columnwidth]{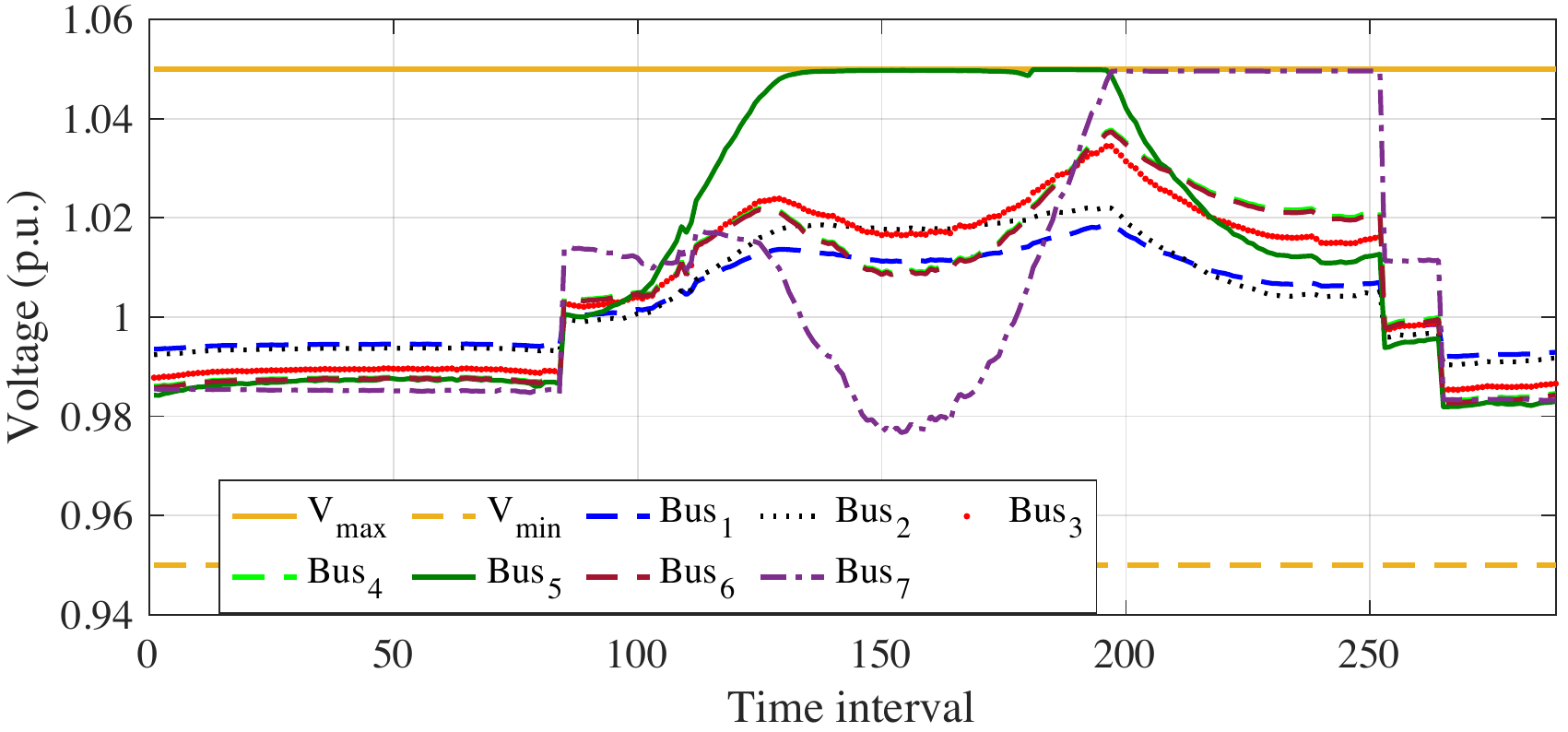}
\caption{Bus voltage profiles of the ETS 3.}
\label{fig:voltages without yn(t)}
\end{figure}

Fig.~\ref{fig:voltages baseline} illustrates the bus voltages of the baseline system.
As shown, due to excess PV power generation during the day, all buses of the feeder experience over-voltage conditions and in the afternoon, due to peak energy demand, the lower voltage limit is breached at buses 3-7. 
Figs.~\ref{fig:voltages all paths}, \ref{fig:voltages without ln(t)} and \ref{fig:voltages without yn(t)} depict that, by exploiting the CES system charging-discharging with PV power generation, all three ETSs can regulate the voltage excursions in the DN that occurs in the baseline system.

\section{Conclusions} \label{conclusion}
We have studied the extent to which the integration of residential photovoltaic (PV) power with a community energy storage (CES) system can minimize distribution network (DN) power loss and energy costs, by comparing three different energy trading systems (ETSs). A network-constrained multi-objective optimization framework has been developed to study the trade-offs between the energy cost and DN power loss reductions while satisfying the DN voltage and current flow limits. Numerical results highlighted that energy trading between the PV users and the CES system can create a better trade-off between the network power loss and energy cost reductions while reducing peak electricity demand on the grid. Future work includes exploring methods for optimal sizing and location of the CES system and combining the ETS operations with behind-the-meter energy storage systems.  
 



\end{document}